\documentclass[%
 twocolumn,aps,prl,
 amsmath,amssymb
]{revtex4-1}

\usepackage{graphicx}
\usepackage{dcolumn}
\usepackage{bm}

\begin{document}

\preprint{APS/123-QED}

\title{Limits on reconstruction of dynamical networks}

\author{Jiajing Guan, Tyrus Berry, Timothy Sauer}
\affiliation{%
 George Mason University\\Fairfax, VA 22030
}%
%

\begin{abstract}
An observability condition number is defined for physical systems modeled by network dynamics. Assuming the dynamical equations of the network are known and a noisy trajectory is observed at a subset of the nodes, we calculate the expected distance to the nearest correct trajectory as a function of the observation noise level, and discuss how it varies over the unobserved nodes of the network. When the condition number is sufficiently large, reconstructing the trajectory from observations from the subset will be infeasible. This knowledge can be used to choose an optimal subset from which to observe a network.
\end{abstract}

\maketitle



The study of network dynamics is increasingly useful in modeling physical processes. Networks present a fascinating departure from generic dynamical systems due to the constraints imposed on direct communication between nodes, resulting in complicated dynamics and nontrivial bifurcation structures \cite{golubitsky2006nonlinear,sauer2004reconstruction,boccaletti2006complex,newman2011structure}.  Modeling by networks has become an important topic in almost every area of physical and biological science, including distributed mechanical processes, weather and climate, and metabolic, genomic and neural networks.

A crucial aspect of studying distributed systems is the difficulty of finding generic observables that facilitate reconstruction of the entire collective dynamics of the network. The theory of observability was pioneered for linear dynamics by Kalman \cite{kalman1959general}. For nonlinear dynamics, the theory of attractor reconstruction \cite{Takens,SYC} provides hope that for generic observables of sufficiently high dimension, the dynamics can be reconstructed. Although observations at single or even multiple nodes of a network may not be provably generic, the results of Joly  \cite{joly2012observation} show that some aspects of observability may  be present by observing even a single node in a strongly connected network, i.e. a network for which every node is downstream from every other node. 
Observability in both linear and nonlinear networks is a topic of intense  recent interest  \cite{letellier2005relation,letellier2005graphical,letellier2009symbolic,liu2013observability,sendina2016observability,whalen2015observability} and has close connections to controllability  \cite{lin1974structural,liu2011controllability,wang2014controllability}. 

However, observability in theory does not guarantee a satisfactory reconstruction from data collected from a sparsely or weakly--connected network, or far from target nodes, even in the  case where the equations of motion are known. To date, even in this more tractable scenario, surprisingly little in the way of general practical requirements  have been developed for inferring information from measurements. A critical obstruction is the presence of noise in the observations, and how this is magnified in efforts to reconstruct the dynamics. In this article, we offer a definition of observability condition number for network reconstruction, show that it has good asymptotic properties for limiting cases such as full observability, and exhibit its behavior for some relevant examples. The main conclusion is that for practical use of network reconstruction techniques, theoretical observability may be only a first step, and that a condition number measuring error magnification may fundamentally govern the limits of reconstructibility.

In the following, we denote by $S$ a subset of observing nodes or variables of a dynamical network. Let $X$ be a network node whose trajectory needs to be reconstructed. We consider an ergodic trajectory of a compact attractor which is observed with noise, and consider the trajectory reconstruction error at one node $X$ of the network. 

In this scenario, we conjecture that there is a constant $\kappa$ depending on $S,X$, and the dynamics, such that in the low noise limit, the expected error of reconstructing a length-$N$ trajectory satisfies
\begin{equation}
{ \mathbb E} \left[  \frac{\rm reconstruction\ error\  per\  step\ at\ X}{\rm observation\  error\  per\  step\ at\ S}
\right] \sim \frac{\kappa}{\sqrt{N}}.
\end{equation}
for some constant $\kappa$. More precisely, we claim that for trajectories observed with noise level $\sigma$, the limit
\begin{equation}
\lim_{N\to\infty, \sigma\to 0}{ \sqrt{N}\ \mathbb E} \left[  \frac{{\rm RMS\ \ error\ \  per\ \ step\ at} \ X}{{\rm obs\ \ error\ \ per\ \ step\ at}\ S}
\right] = \kappa,
\end{equation}
exists. Equivalently, $\kappa$ can be defined using the Euclidean 2-norm
\begin{equation} \label{e3}
\kappa=\lim_{N\to\infty, \sigma\to 0}{ \mathbb E} \left[  \frac{||e||_2}{\sigma}
\right],
\end{equation}
where $e=\{e_1,\ldots, e_N\}$ denotes the trajectory error at $X$ with component $e_i=z_i-x_i$ in terms of the exact trajectory $x_i$ and the reconstructed trajectory $z_i$. We call $\kappa = \kappa_{S,X}$ the {\it observability condition number} of node $X$ observed by $S$. This is a single constant that encapsulates the ability to reconstruct the dynamics at $X$ from the subset $S$.

The study of condition number as a measure of controllability and observability is classical, beginning with Friedland \cite{friedland1975controllability} in the context of linear systems. Here we consider the nonlinear case, and append the asymptotics expressed in (\ref{e3}). In addition,  we describe a direct computational approach to approximating the condition number: Small observational noise is added to a length $N$ trajectory of the dynamical system and a variational data assimilation technique is used to reconstruct the nearest exact trajectory. The ratio of output (reconstruction) error to input (observation) error is $\kappa_{S,X}/\sqrt{N}$. 

The fundamental importance of the existence of a universal quantity $\kappa_{S,X}$, independent of trajectory length, is that it allows us to compare the effects of various observation sets $S$ at node $X$. This has direct implications for sensor placement in general systems, such as the positioning of weather buoys or the location of electrodes in a neural assembly.

\begin{figure}
 \begin{center}
{\includegraphics[width=.4\linewidth]{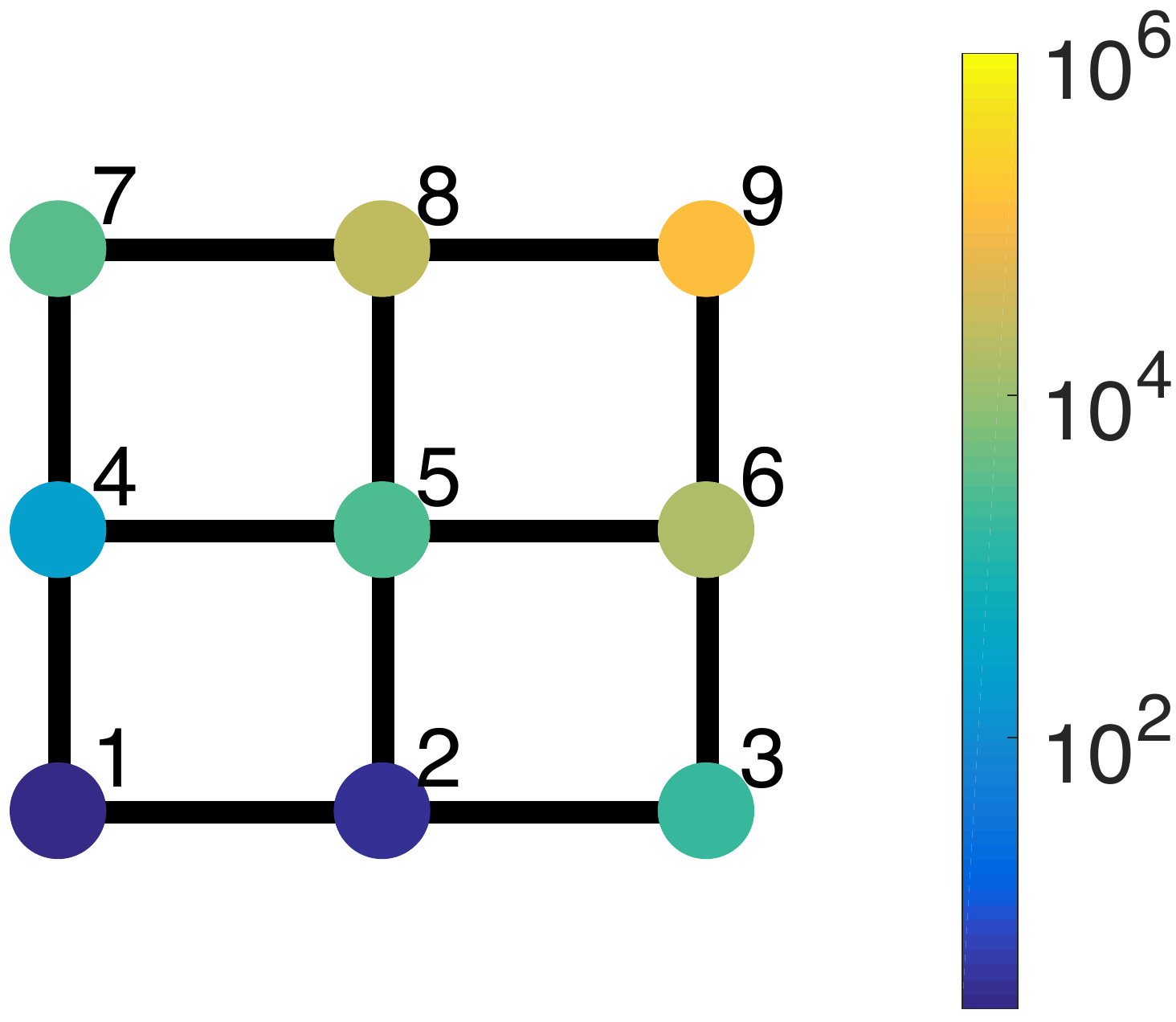}}\hspace{.02\linewidth}
{\includegraphics[width=.56\linewidth]{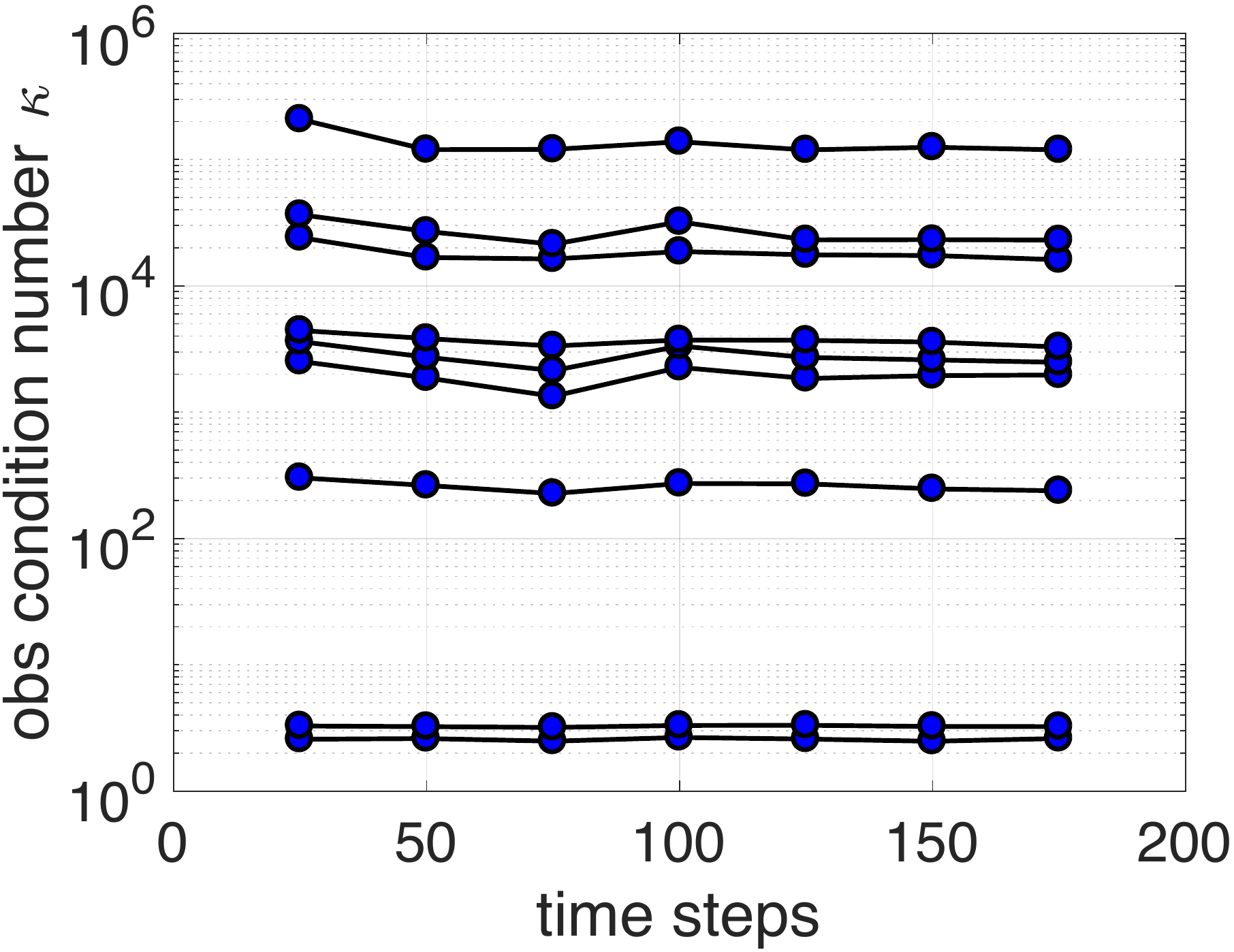}}
\hspace*{-.0\linewidth} (a) \hspace*{.4\linewidth} (b)\\
{\includegraphics[width=.28\linewidth]{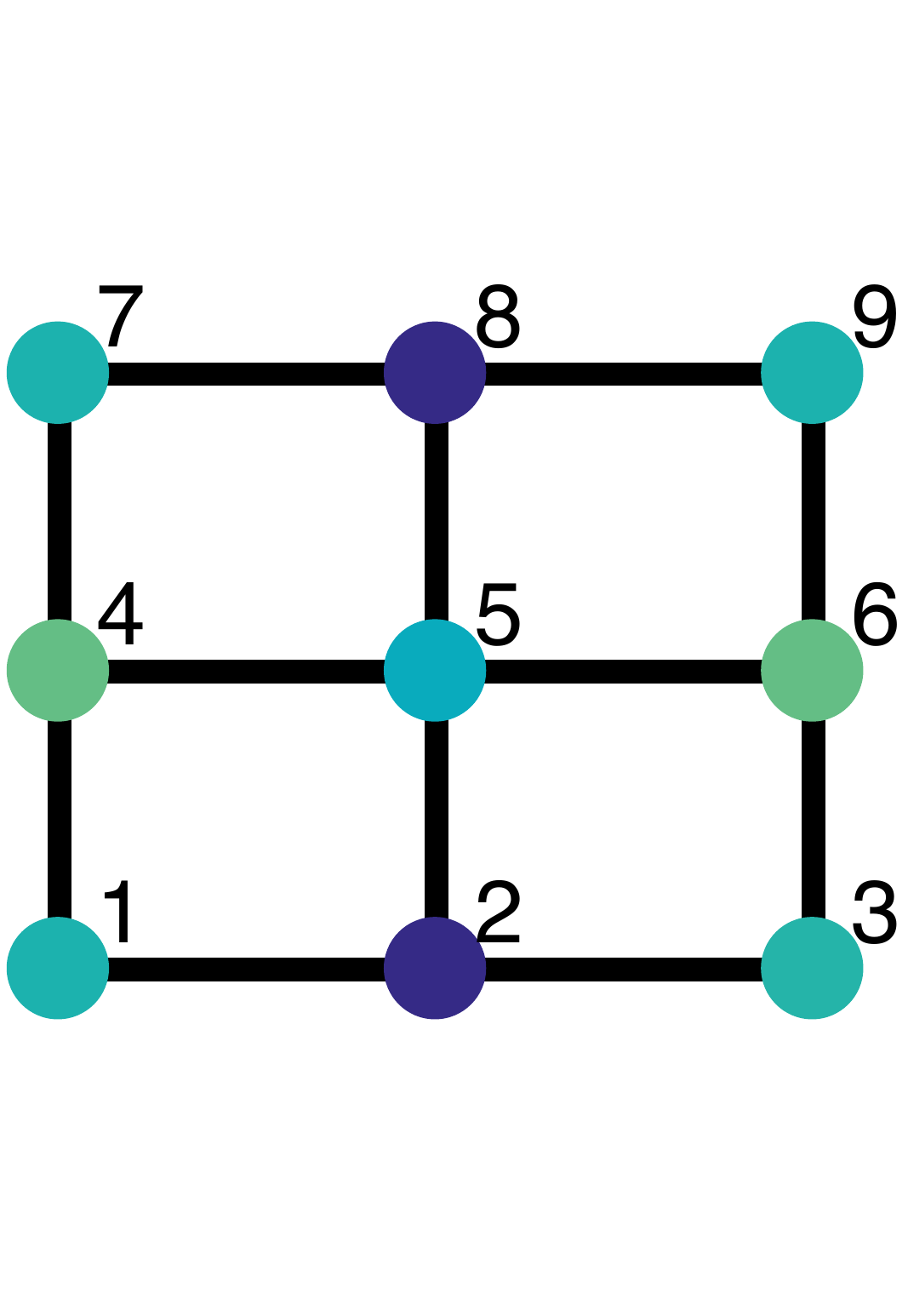}}\hspace{.13\linewidth}
{\includegraphics[width=.56\linewidth]{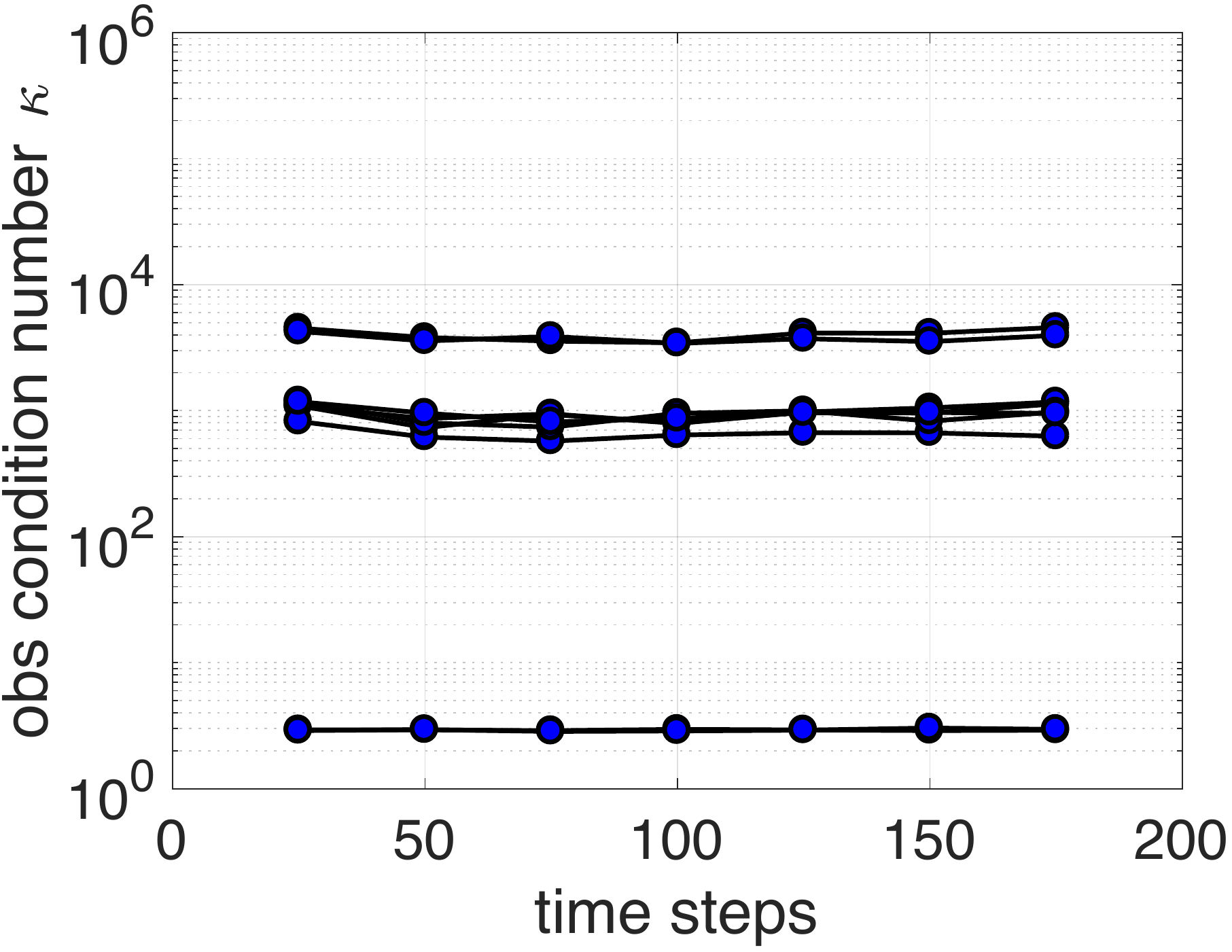}}
\hspace*{-.0\linewidth} (c) \hspace*{.4\linewidth} (d)\\
\end{center}
\caption{Undirected network of nine nodes with dynamics (\ref{e5}) and adjacency matrix $A$ as shown. (a) Node $X$ is shaded (color online) according to $\kappa_{S,X}$ where $S=\{1,2\}$. (b) Estimates of $\kappa_{S,X}$ in (a) as a function of trajectory length.(c) Same as (a), but $S=\{2,8\}$.  (d) Estimates of $\kappa_{S,X}$ in (c). } \label{f1}
\end{figure}

We begin by establishing the formula (\ref{e3}) for the case of completely observed linear dynamics, where $\kappa$ exists and equals $1$. By ``completely observed'', we mean that the subset $S$ of observed variables includes all variables. Consider first the scalar case and assume the dynamics $f(x)=ax$. Let $\{x_1,\ldots,x_N\}$ be a trajectory under $f$,  so that $x_i=f^{(i-1)}(x_1)$. Let $y_i=x_i+\epsilon_i$ be the (completely observed) trajectory observed with i.i.d. observation noise $\epsilon_i$ of mean $0$ and covariance $\Sigma(\epsilon) = \sigma^2{\mathbb I}_{N\times N}$ for some $\sigma>0$. Consider the exact trajectory $\{z_1,\ldots,z_N\}$ where $z_i=f^{(i-1)}(z_1)$, that minimizes the sum squared error 
\begin{equation} \label{e7}\sum_{i=1}^N (z_i-y_i)^2 = \sum_{i=1}^N (a^{i-1}z_1-a^{i-1}x_i-\epsilon_i)^2.\end{equation}
In the sense of least squares, the $\{z_i\}$ trajectory is the one closest to the observations.
Setting the derivative with respect to $z_1$ to zero and solving yields
\begin{eqnarray} \label{e4}
0 &=& z_1\sum_{i=1}^N a^{2(i-1)}- x_1\sum_{i=1}^N a^{2(i-1)}-\sum_{i=1}^N a^{i-1}\epsilon_i\nonumber\\
z_1 &=& x_1+\frac{\sum_{i=1}^N a^{i-1}\epsilon_i}{\sum_{i=1}^N a^{2(i-1)}}.
\end{eqnarray}
The square of the numerator of (\ref{e3}) is the expected squared error of the reconstructed trajectory
 $\{z_1,\ldots,z_N\}$ compared with the original trajectory  $\{x_1,\ldots,x_N\}$, or using $z_i-x_i = a^{i-1}(z_1-x_1)$ and (\ref{e4}),
\begin{eqnarray*}
{\mathbb E} \left[
\sum_{i=1}^N (z_i-x_i)^2\right] &=&{\mathbb E} (z_1-x_1)^2  \sum_{i=1}^N a^{2(i-1)} \\
&=& {\mathbb E} \left[ \left( \frac{ \sum_{i=1}^N a^{i-1}\epsilon_i}{ \sum_{i=1}^N a^{2(i-1)} }\right)^2\right]  \sum_{i=1}^N a^{2(i-1)}\\
&=& \frac{\sum_{i=1}^N a^{2(i-1)}{\mathbb E}(\epsilon_i^2)}{\sum_{i=1}^N a^{2(i-1)}} = \sigma^2
\end{eqnarray*}
where we used the fact that the noises $\epsilon_i$ are uncorrelated. Dividing by the observation noise level $\sigma$, we conclude that $\kappa = 1$ for the completely observed case. 

The scalar case can be extended to linear dynamics $f(x)=Ax$ for a symmetric matrix by diagonalizing $A$ and applying the above argument componentwise. Furthermore, an argument based on the multiplicative ergodic theorem \cite{MET} allows us to extend the fact that $\kappa=1$ for the completely observed case to general nonlinear dynamics. 

  \begin{figure}
 \begin{center}
{\includegraphics[width=.14\linewidth]{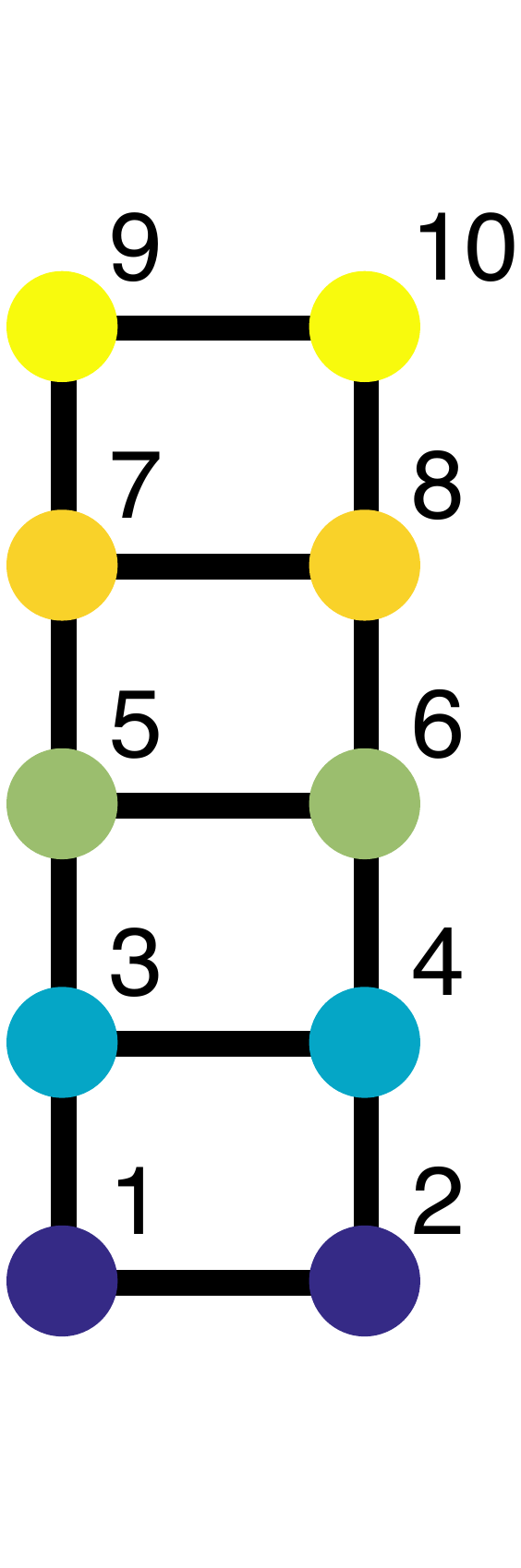}}\hspace{.03\linewidth}
{\includegraphics[width=.6\linewidth]{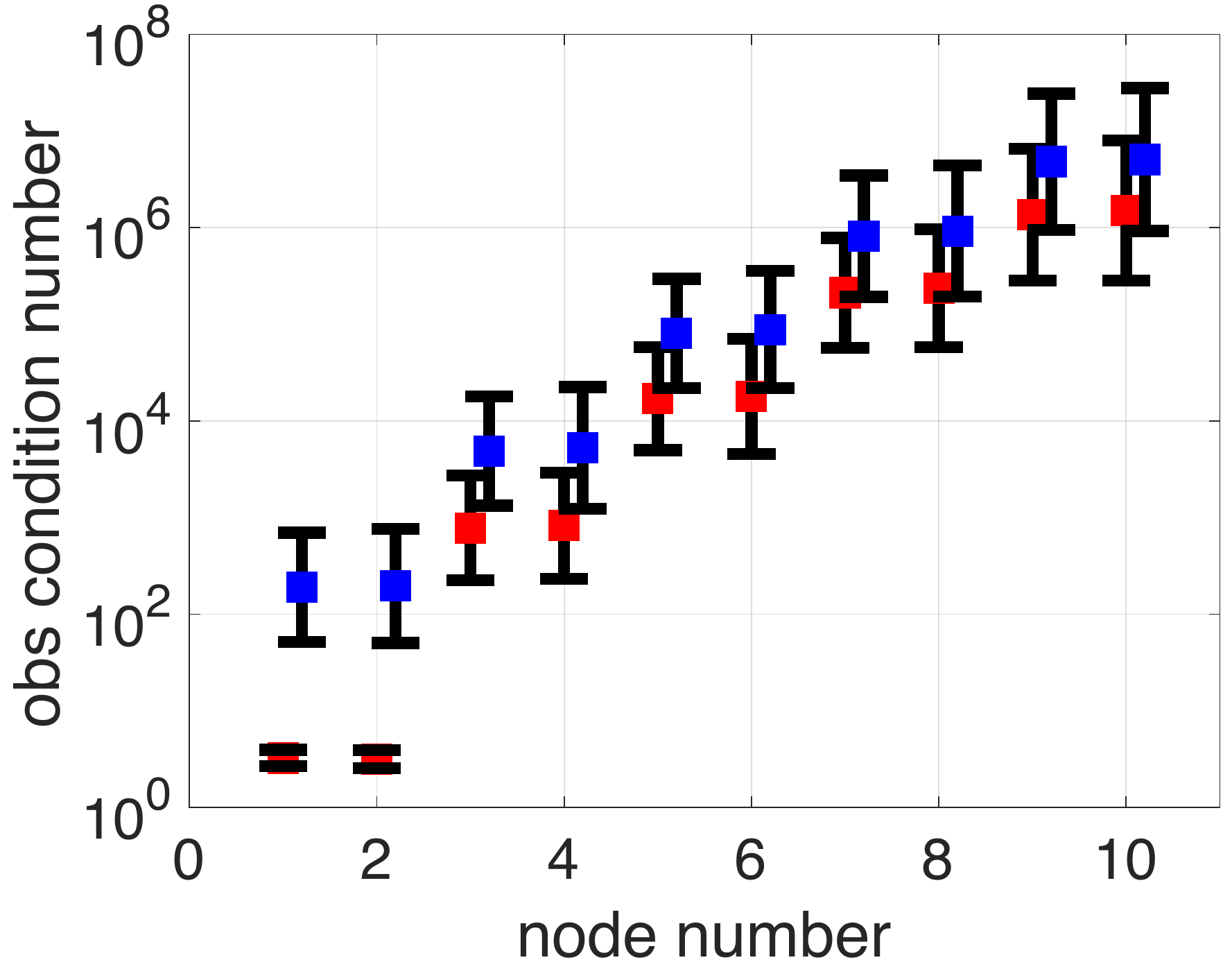}}\\
\hspace*{-.0\linewidth} (a) \hspace*{.4\linewidth} (b)\\
\end{center}
\caption{(a) A network of 10 nodes with dynamics from (\ref{e5}). (b) Plot of $\kappa_{S,X}$ for nodes 1 -- 10 when observed at $x$-coordinates of subset $S=\{1,2\}$. The means and standard deviations for the $x^j$ (left) and $y^j$ (right) variables are shown.}\label{f2}
\end{figure}

For partial observations, such as observing at a proper subset $S$ of nodes of a network, $\kappa_{S,X}$ will be substantially greater than $1$, which is the focus of this article. 
As an illustrative example, consider the undirected network of nine nodes illustrated in Fig. \ref{f1}, where the update equations at node $j$ follow the nonlinear discrete dynamics
\begin{equation} \label{e5}
x^j_{i+1} = a_j\cos x^j_i + b_jy^j_i + c\sum_{k=1}^9 A_{jk}x^k_i,\ \ \ \ 
y^j_{i+1} = x^j_i
\end{equation}
where $a_j=2.2, b_j=0.4$ for $j=1,\ldots,9$ and $A=\{A_{jk}\}$ is the (symmetric) adjacency matrix of the network. The discrete dynamical map used here at each node is a variant of the classical Henon map \cite{henon1976two} that is suitable for distributed dynamics.

We describe an algorithm to compute $\kappa$ as in (\ref{e3}) from a general network. Generate an exact trajectory $\{x_1,\ldots, x_N\}$, which is observed by a function $h(x)$ plus Gaussian observational noise with variance $\sigma^2$ at each point of the trajectory to get $\{y_1,\ldots, y_N\}$, where $y_i = h(x_i)+\epsilon_i$. In the examples to follow, $h$ will represent observing at the subset $S$ of nodes. We apply a variational data assimilation method to the inexact $y_i$ observations to find an exact trajectory $\{z_1,\ldots, z_N\}$ that minimizes the least squares difference between the $y_i$ and $z_i$ trajectories, analogous to (\ref{e7}). To accomplish this, we applied a Gauss-Newton iteration that enforces the exactness of the $z_i$ trajectory while minimizing the observation difference. More precisely, we minimize
\begin{equation} \label{e8} \frac{1}{q^2}\sum_{i=1}^{N-1} (f(z_i)-z_{i+1})^2+\frac{1}{r^2}\sum_{i=1}^{N-1} (h(z_i)-y_i)^2
\end{equation}
where $q$ and $r$ are weights that specify the trajectory noise and observation noise tolerances, respectively. We use $q << r$ to ensure that the $z_i$ trajectory is effectively exact, at least relative to the observation errors. At the conclusion of the Gauss-Newton iteration for $z_i$, we compute the errors $e_i=z_i-x_i$ and the approximation (\ref{e3}) to $\kappa_{S,X}$. 

The results of this algorithm applied to the network in Fig.~\ref{f1}(a), observed with the $x$-coordinates at the set $S=\{1,2\}$, are shown in Fig.~\ref{f1}(b). The nine traces correspond to each of the nine network nodes $X$. The two observed nodes are at the bottom, and remaining traces show various levels of $\kappa_{S,X}$. In this example, the asymptotic $N\to\infty$ limit in (\ref{e3}) is reached for relatively short trajectory lengths. The  nodes in Fig.~\ref{f1}(a) are colored according to the respective observability condition numbers.

The fact that $\kappa_{S,X}$ can be arbitrarily large is illustrated by undirected networks such as Fig.~\ref{f2}. The equations (\ref{e5}) are used, and the observing set is $S=\{1, 2\}$. As expected, the resulting condition numbers grow with the distance from the observing set. However, calculating large $\kappa_{S,X}$ is delicate, as we discuss below.

\begin{figure}
 \begin{center}
{\includegraphics[width=.55\linewidth]{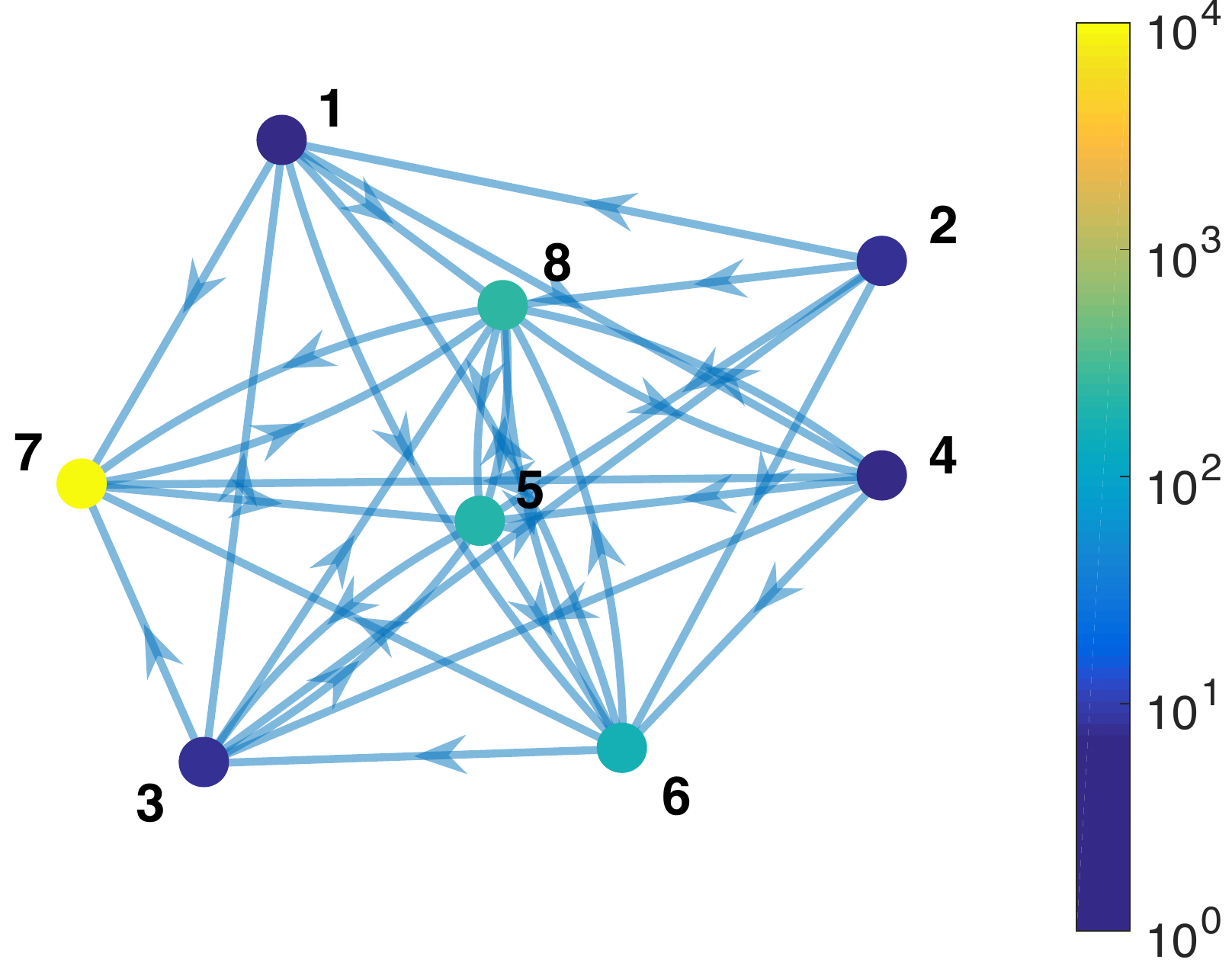}}\hspace{.1\linewidth}\\
(a) \\
{\includegraphics[width=.62\linewidth]{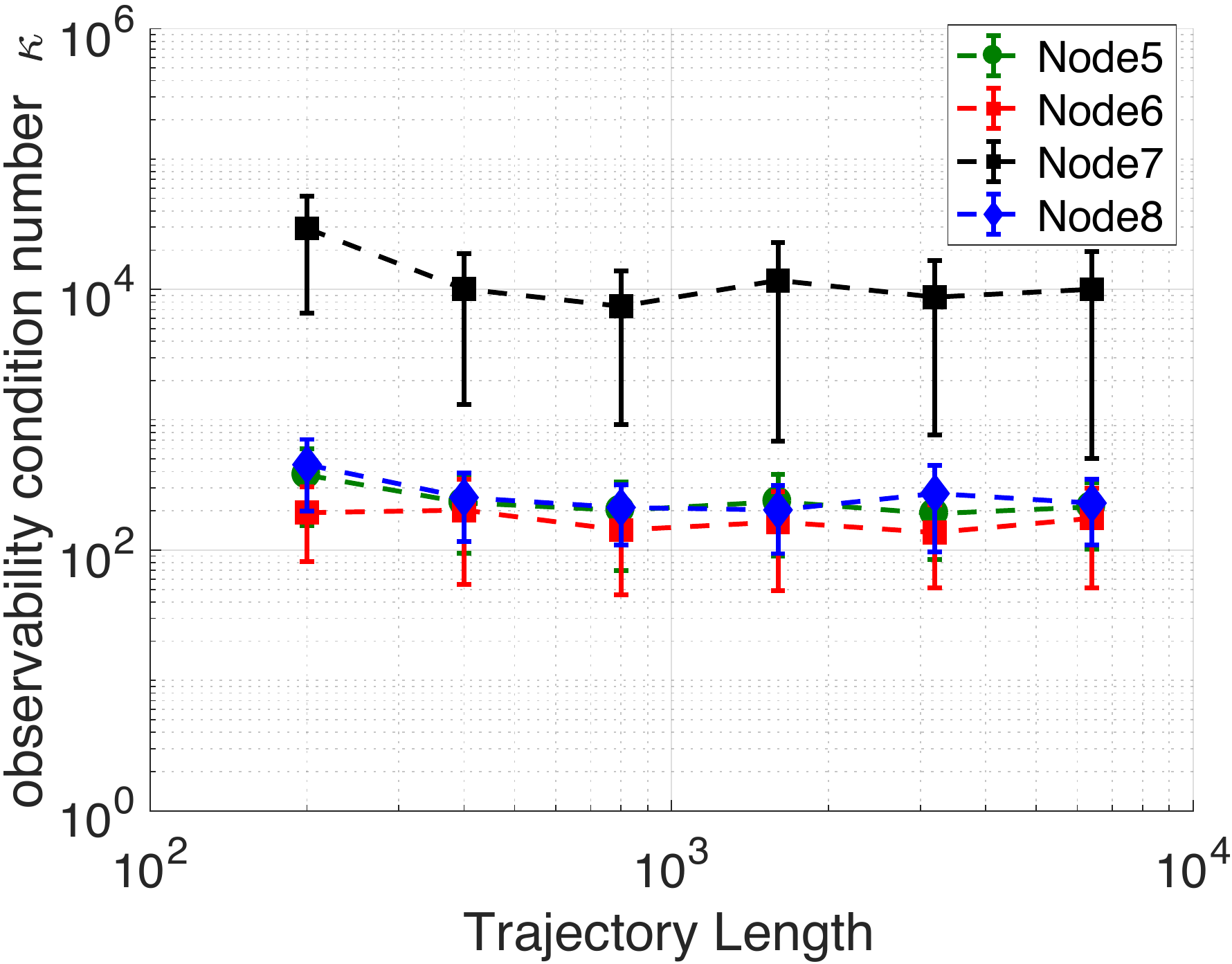}}\\
(b)
\end{center}
\caption{(a) $\kappa_{S,X}$ for Fitzhugh-Nagumo network. (b) Estimates of $\kappa_{S,X}$ versus trajectory length. Although not obvious from the adjacency matrix, node $7$ is more difficult to observe from the set $S=\{1, 2, 3, 4\}$ than the remaining nodes.} \label{f3}
\end{figure}

As an example of a discretely-sampled continuous dynamical sysfem we built a directed network of Fitzhugh-Nagumo neurons \cite{fitzhugh1961impulses,nagumo1962active} as shown in Fig.~\ref{f3}(a)
\begin{eqnarray} \label{e6}
v_j' &=& -w_j+dv_j-v_j^3/3+I+g\sum_{k=1}^9 A_{jk}v_k\nonumber\\
w_j' &=& a-bw_j+cv_j
\end{eqnarray}
where the parameters were varied by about $5\%$ from $a=0.42, b=0.8, c=0.08, d=0.01$ and $I=-0.025$ among the nodes. The system was observed at nodes 1 -- 4 at a step size $\Delta t = 0.1$, and small observation noise was added at each step. Fig.~\ref{f3}(b) shows the observability condition number calculated at the remaining four nodes.

\begin{figure}
 \begin{center}
{\includegraphics[width=.85\linewidth]{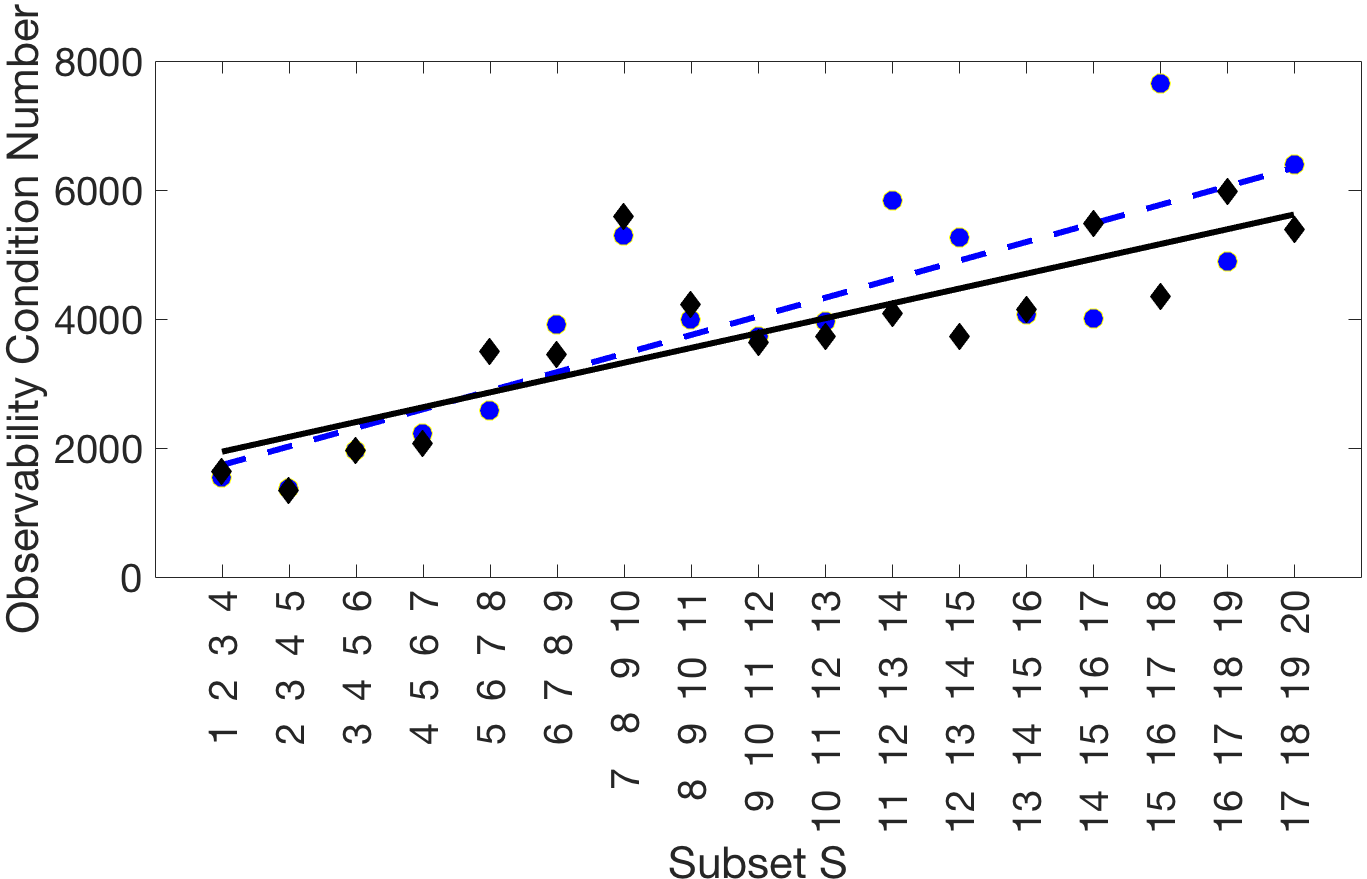}}\\
(a)\\
{\includegraphics[width=.85\linewidth]{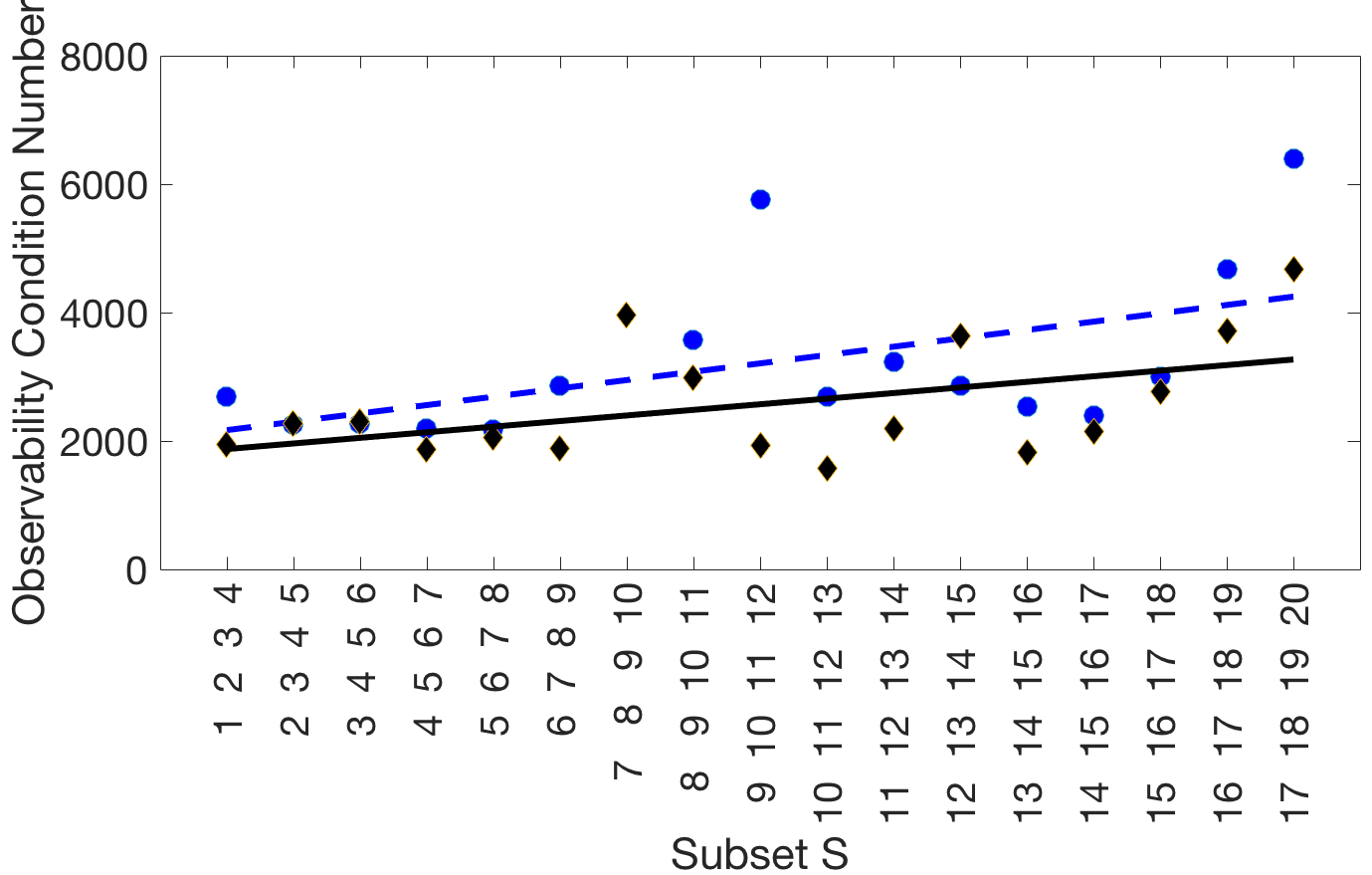}}\\
(b)
\end{center}
\caption{(a) Mean observability condition number of scale-free network, as a function of four-member observation subsets $S$. Round markers (dashed lines) denote nodes sorted by descending degree; diamond markers (solid lines) denote nodes sorted by descending closeness centrality. (b) Same for Erd\"os-Renyi network.} \label{f4}
\end{figure}

We tested the computation of $\kappa_{S,X}$ in two other networks to compare the effects of hubness. The example in Fig.~\ref{f4}(a) is a scale-free network of 20 nodes. After sorting the nodes in descending order by degree (resp., centrality),  $\kappa_{S,X}$ was computed for $S$ equal to the first four, second four, etc. The more sparsely-connected observer sets lead to much increased mean $\kappa_{S,X}$. The same analysis, but for an Erd\"os-R\'enyi network of the same size, is done in Fig.~\ref{f4}(b). Note that values of $\kappa_{S,X}$ vary much less with the choice of the subset $S$.

To conclude, we note that special care must be taken to carry out the Gauss-Newton iteration $$z^{k+1} = z^k-(DR^TDR)^\dagger (DR)^TR(z^k),$$ where we have set $R=[R_1,\ldots, R_n]$ and $ R_1(z)^2+\ldots+R_n(z)^2$ is the sum (\ref{e8}) to be minimized. Clearly, the accuracy of the solution of the iteration will be dependent on the conditioning of the Gauss-Newton problem, in this case the (traditional) condition number $C$ of the matrix $(DR)^TDR$. Since the errors added are of size $\sigma<<1$, the residuals $R_i$ can be expected to be on the same order. Thus one may run out of correct significant digits if $C/\sigma > \epsilon_{\rm mach}^{-1} \approx 10^{16}$ for double precision computations. 

\begin{figure} 
    \includegraphics[width=.75\linewidth]{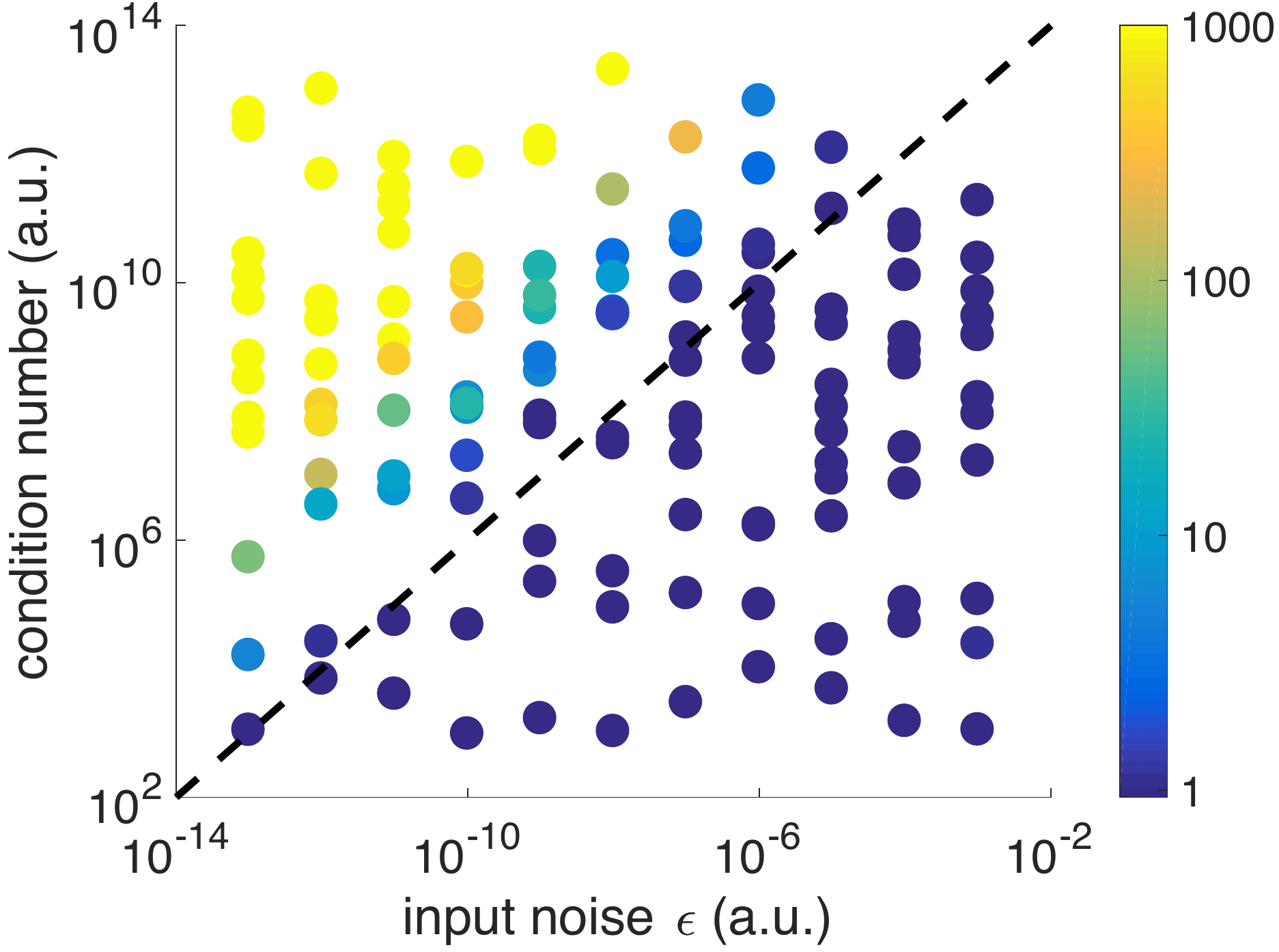}
    \caption{Markers are shaded by estimate of $\kappa_{S,X}$ and plotted versus input noise and the condition number $C$ of the Gauss-Newton iteration, for completely observed network of random $2\times 2$ matrices. Trajectory lengths vary from 20 (lower points) to 240 (upper points). Calculations below the dashed line $C = 10^{16}\sigma$ are reliable. } \label{f5}
\end{figure}

In Fig.~\ref{f5} we explore this issue. For simplicity, we consider a discrete stochastic map that multiplies by a random $2\times 2$ matrix at each step. We observe at both phase variables, so that the system is completely observable and $\kappa_{S,X}=1$.
 Each marker represents a calculation of $\kappa_{S,X}$ where $S$ is both variables, i.e. completely observed, and $X$ is one of the variables. As shown above, in this case $\kappa_{S,X}=1$. The markers correspond to trajectories of length between 20 and 240, going from bottom to top, with input noise $\sigma$. The vertical axis denotes the condition number $C$ of $(DR)^TDR$. The color of the marker corresponds to the observability condition number $\kappa_{S,X}$. The dashed line is drawn at $C/\sigma = 10^{16}$. Note that as the trajectories become longer, $C$ becomes larger and when the dashed line is passed, $\kappa_{S,X}$ is incorrectly determined (in some cases by a factor of more than 1000), due to lack of significant digits caused by ill-conditioning. To avoid this difficulty, the length of trajectories must be limited to the safe area below the dashed line. Alternatively,  computations beyond double precision could be used.

In this article, we have introduced the concept of observability condition number $\kappa_{X,S}$ that has a consistent asymptotic definition in the limit of long ergodic trajectories and the limit of small noise. We have shown that the definition in relatively straightforward to compute in multidimensional systems. This settles a fundamental, long-standing problem in network dynamics, namely where to locate a minimal set of sensors to measure remote dynamics. Computation of $\kappa_{S,X}$ allows a direct comparison of all options.
 In particular, an exhaustive enumeration among subsets $S$ to find the minimum  mean or maximum over the network is feasible for moderate-sized networks, and establishes a guiding principle for large networks where exhaustive search may not be feasible.

\begin{acknowledgments}This work was partially supported by NSF grant DMS-1723175.
\end{acknowledgments}

\bibliography{obscond}

\end{document}